\documentclass[11pt,leqno]{article}
\usepackage{amsmath,amssymb,color}
\usepackage[colorlinks=true,breaklinks=true,linkcolor=blue,citecolor=red]{hyperref}

\font\texteusm eusm10 scaled1095 
\font\scripteusm eusm7
\font\scriptscripteusm eusm5
\newfam\eusmfam
\textfont\eusmfam\texteusm
\scriptfont\eusmfam\scripteusm
\scriptscriptfont\eusmfam\scriptscripteusm

\newtheorem{Definition}{Definition}  
\newtheorem{Lemma}{Lemma}  
\newtheorem{Remark}{Remark}[section]
\newtheorem{Example}{Example}[section]

\newtheorem{Corollary}{Corollary}

\newtheorem{Theorem}{Theorem}  
{\catcode `\@=11 \global\let\AddToReset=\@addtoreset}
\AddToReset{equation}{section}
\AddToReset{Conjecture}{section}
\AddToReset{Theorem}{section}
\AddToReset{Proposition}{section}
\AddToReset{Corollary}{section}
\AddToReset{Lemma}{section}
\AddToReset{Definition}{section}
\AddToReset{figure}{section}
\begin{document}
 
\title{On the absence of global weak solutions for a nonlinear time-fractional  Schr\"odinger equation} 

\author{Munirah Alotaibi$^1$,  Mohamed Jleli$^1$, Maria Alessandra Ragusa$^{2,3}$, Bessem Samet$^1$}

\date{\small{$^1$Department of Mathematics, College of Science, King Saud University,
 P.O. Box 2455, Riyadh, 11451, Saudi Arabia\\
maalotaiby@pun.edu.sa (M. Alotaibi); jleli@ksu.edu.sa (M. Jleli); bsamet@ksu.edu.sa (B. Samet)\\
\vspace{0.2cm}
$^2$Department of Mathematics and Computer Science, University of Catania,
 95124 Catania, Italy\\\vspace{0.2cm}
 $^3$RUDN University, 6 Miklukho, Maklay St, Moscow 117198, Russia\\
 mariaalessandra.ragusa@unict.it (M.A. Ragusa)}}
\maketitle 

\begin{abstract} 
In this paper, an initial value problem for a nonlinear time-fractional  Schr\"odinger equation  with a singular logarithmic potential term is investigated. The considered problem involves  the left/forward Hadamard-Caputo fractional derivative with respect to the time variable.  Using the test function method  with a judicious choice of the test function, we obtain sufficient criteria for the absence of global weak solutions. 
 \end{abstract}

\section{Introduction}\label{sec1}

This paper is concerned with the nonexistence of global weak solutions to the initial value problem for the time-fractional  Schr\"odinger equation 
\begin{eqnarray}\label{P}
\left\{\begin{array}{llll}
i^\alpha \mathcal{D}_{a|t}^\alpha u+\Delta u	=\lambda \displaystyle \left(\ln \frac{t}{a}\right)^\gamma  |u|^p,\quad t>a,\,\, x\in \mathbb{R}^N,\\ \\
u(a,x)=f(x),\quad x\in \mathbb{R}^N,
\end{array}
\right.
\end{eqnarray}
where $u=u(t,x)$ is the complex-valued unknown function, $N\geq 1$, $a>0$, $i\in \mathbb{C}$ is  the imaginary unit ($i^2=-1$), $0<\alpha<1$, $i^\alpha= e^{i\alpha \pi/2}$,  $\mathcal{D}_{a|t}^\alpha$ is the left/forward Hadamard-Caputo fractional derivative of order $\alpha$ with respect to the time variable $t$ (see Section \ref{sec2}), $\Delta$ is the Laplacian operator with respect to the space variable $x$, $\lambda\in \mathbb{C}\backslash\{0\}$,  $\gamma\in \mathbb{R}$, $p>1$, and $f\in L^1(\mathbb{R}^N,\mathbb{C})$. Namely, we are interested in obtaining sufficient conditions under which \eqref{P} admits no global weak solution, in the sense that will be specified later. 

The theory of fractional calculus has received a great attention from several researchers 
working in various disciplines. Namely, it was shown that many real-world phenomena can be better modeled using fractional operators, see e.g. \cite{AB,BAG,HI,IO,LAS,MA,MAI}, and the references therein. Due to this fact, the study of fractional evolution equations   has become increasingly popular.  In particular, the study of time-fractional schr\"odinger equations in both, theoretical and numerical aspects, has been attracted a great deal of attention, see e.g. \cite{AR,DO,GU,LAS2,LI,NAB,RA,WA,WA2}, and the references therein. 

For any complex number $z\in \mathbb{C}$, we denote by $\Re z$ and $\Im z$  the real and the imaginary parts of $z$, respectively.

In \cite{KI-NA}, Kirane and Nabti considered the initial value problem for the nonlocal in time nonlinear Schr\"odinger equation
\begin{eqnarray}\label{PKI-NA}
\left\{\begin{array}{llll}
i\displaystyle\frac{\partial u}{\partial t}+\Delta u	=\displaystyle\frac{\lambda}{\Gamma(\alpha)} \int_0^t (t-s)^{\alpha-1} |u(s,x)|^p\,ds,\quad t>0,\,\, x\in \mathbb{R}^N,\\ \\
u(0,x)=f(x),\quad x\in \mathbb{R}^N,
\end{array}
\right.
\end{eqnarray}
where $0<\alpha<1$   and $\Gamma(\cdot)$ is the Gamma function.  It was shown that, if 
$$
\Re f\in L^1(\mathbb{R}^N,\mathbb{R}),\,\, \Im \lambda\int_{\mathbb{R}^N} \Re f(x)\,dx>0\quad\mbox{or}\quad \Im f\in L^1(\mathbb{R}^N,\mathbb{R}),\,\, \Re \lambda\int_{\mathbb{R}^N} \Im f(x)\,dx<0,
$$
and 
$$
1<p\leq 1+\frac{2(\alpha+1)}{N-2\alpha},\quad  N>2\alpha,
$$
then \eqref{PKI-NA} has has no global weak solution. 

In \cite{ZH},  Zhang et al. considered the nonlinear time-fractional Schr\"odinger equation 
\begin{eqnarray}\label{P-Zhang}
\left\{\begin{array}{llll}
i^\alpha \,{}^C\!D_{0|t}^\alpha u+\Delta u	=\lambda  |u|^p,\quad t>0,\,\, x\in \mathbb{R}^N,\\ \\
u(0,x)=f(x),\quad x\in \mathbb{R}^N,
\end{array}
\right.
\end{eqnarray}
where $0<\alpha<1$ and  ${}^C\!D_{0|t}^\alpha$ is the Caputo fractional derivative of order $\alpha$ with respect to the time variable $t$  (see \cite{KIL} for the definition of Caputo fractional derivative). It was shown that, if $1<p<1+\frac{2}{N}$,   $f\in L^1(\mathbb{R}^N,\mathbb{C})$, and 
$$
\Re \lambda\int_{\mathbb{R}^N}F_1(x)\,dx>0\quad\mbox{or}\quad \Im \lambda\int_{\mathbb{R}^N}F_2(x)\,dx>0,
$$
where
$$
F_1(x)=\cos\left(\frac{\pi\alpha}{2}\right)\Re f(x)-\sin\left(\frac{\pi\alpha}{2}\right) \Im f(x)
$$
and
$$
F_2(x)=\cos\left(\frac{\pi\alpha}{2}\right)\Im f(x)+\sin\left(\frac{\pi\alpha}{2}\right) 
\Re f(x),
$$
then \eqref{P-Zhang} admits no global weak solution. Let us mention that $1+\frac{2}{N}$ is the Fujita critical exponent for the semilinear heat equation  $\frac{\partial u}{\partial t}-\Delta u=|u|^p$, $t>0$, $x\in \mathbb{R}^N$ (see Fujita \cite{Fujita}).

Let us mention that in the limit case $\alpha\to 1^-$, problem \eqref{P-Zhang} reduces to the nonlinear time Schr\"odinger equation (see e.g. \cite{KIL}) 
\begin{eqnarray}\label{P-alpha1}
\left\{\begin{array}{llll}
i \displaystyle\frac{\partial u}{\partial t}+\Delta u=\lambda  |u|^p,\quad t>0,\,\, x\in \mathbb{R}^N,\\ \\
u(0,x)=f(x),\quad x\in \mathbb{R}^N.
\end{array}
\right.
\end{eqnarray}
In \cite{IK}, Ikeda and Wakasugi considered the global behavior of solutions to problem \eqref{P-alpha1}, then they established a finite-time blow-up result of an $L^2$-solution, whenever 
$p\in (1,1+2/N)$. Later, the same problem was  discussed by Ikeda and Inui \cite{IKI}, where they established a small data blow-up result of $H^1$-solution, whenever $p\in (1,1+4/N)$.

In  the above mentioned papers, the time fractional derivative was considered in the Caputo sense.  In this paper, we investigate the nonlinear time-fractional  Schr\"odinger equation \eqref{P}, which involves the left/forward Hadamard-Caputo fractional derivative introduced in \cite{AG}. This fractional differential operator differs from the preceding ones in the sense that the kernel of the integral in its definition contains a logarithmic function.  

The rest of the paper is organized as follows. In Section \ref{sec2}, ...

\section{Preliminaries}\label{sec2}

Let $(a,T)\in \mathbb{R}^2$ be such that $0<a<T$. We denote by $AC([a,T],\mathbb{R})$ the space of  real-valued  absolutely continuous functions on $[a,T]$.

Let $f\in L^1([a,T],\mathbb{R})$. The left-sided and right-sided Riemann-Liouville
fractional integrals of of order $\sigma>0$ of $f$, are defined respectively by (see \cite{KIL})
$$
(I_a^\sigma f)(t)=\frac{1}{\Gamma(\sigma)} \int_a^t (t-s)^{\sigma-1}f(s)\,ds	
$$
and
$$
(I_T^\sigma f)(t)=\frac{1}{\Gamma(\sigma)} \int_t^T (s-t)^{\sigma-1}f(s)\,ds	,
$$
for almost everywhere $t\in [a,T]$.

We have the following integration by parts rule.
\begin{Lemma}[see \cite{KIL}]\label{L1}
Let $\sigma>0$, $q,r\geq 1$, and $\frac{1}{q}+\frac{1}{r}\leq 1+\sigma$ ($q=1$, $r=1$, in the case $\frac{1}{q}+\frac{1}{r}=1+\sigma$).	If $f\in L^q([a,T],\mathbb{R})$ and $g\in L^r([a,T],\mathbb{R})$, then 
$$
\int_a^T (I_a^\sigma f)(t)g(t)\,dt=\int_a^T f(t) (I_T^\sigma g)(t)\,dt.
$$
\end{Lemma}

For $r\geq 1$, we denote by $L^r\left([a,T],\mathbb{R},\frac{1}{t}\,dt\right)$ the weighted Lebesgue space of real-valued measurable functions $f: [a,T]\to \mathbb{R}$ satisfying
$$
\int_a^T |f(t)|^r \frac{1}{t}\,dt<\infty. 
$$
Let $f\in L^1\left([a,T],\mathbb{R},\frac{1}{t}\,dt\right)$.  The left-sided  and right-sided Hadamard fractional integrals of  order $\sigma>0$ of  $f$, are defined respectively by (see \cite{KIL})
$$
(J_a^\sigma f)(t)=\frac{1}{\Gamma(\sigma)} \int_a^t  \left(\ln \frac{t}{s}\right)^{\sigma-1} f(s) \frac{1}{s}\,ds
$$
and 
$$
(J_T^\sigma f)(t)=\frac{1}{\Gamma(\sigma)} \int_t^T  \left(\ln \frac{s}{t}\right)^{\sigma-1} f(s) \frac{1}{s}\,ds,
$$
for almost everywhere $t\in [a,T]$.

\begin{Remark}\label{RK-PR}
It can be easily seen that, if $f\in C([a,T],\mathbb{R})$ and $\sigma>0$, then 
$$
\lim_{t\to a^+}(J_a^\sigma f)(t)=\lim_{t\to T^-}  (J_T^\sigma f)(t)=0.
$$	
\end{Remark}

The following integration by parts rule holds.
\begin{Lemma}\label{L2}
Let $\sigma>0$, $q,r\geq 1$, and $\frac{1}{q}+\frac{1}{r}\leq 1+\sigma$ ($q=1$, $r=1$, in the case $\frac{1}{q}+\frac{1}{r}=1+\sigma$).	If $f\in L^q\left([a,T],\mathbb{R},\frac{1}{t}\,dt\right)$ and $g\in L^r\left([a,T],\frac{1}{t}\,dt\right)$, then 
$$
\int_a^T (J_a^{\sigma} f)(t)g(t)\frac{1}{t}\,dt=\int_a^T f(t) (J_T^{\sigma} g)(t)\frac{1}{t}\,dt.
$$
\end{Lemma}

\noindent{\it Proof.}  Using the change of variable $\tau=\ln s$, we obtain
$$
(J_a^{\sigma} f)(t)=\frac{1}{\Gamma(\sigma)} \int_{\ln a}^{\ln t} (\ln t-\tau)^{\sigma-1}(f\circ \exp)(\tau)\, d\tau,
$$
that is,
\begin{equation}\label{in1}
(J_a^{\sigma} f)(t)=\left(I_{\ln a}^\sigma\,\, f\circ \exp\right) (\ln t).	
\end{equation}
Using the same change of variable, we obtain
\begin{equation}\label{in2}
(J_T^{\sigma}g)(t) =\left(I_{\ln T}^\sigma \,\,g\circ \exp\right) (\ln t).
\end{equation}
By \eqref{in1}, there holds
$$
\int_a^T (J_a^{\sigma} f)(t)g(t)\frac{1}{t}\,dt=\int_a^T  \left(I_{\ln a}^\sigma\,\, f\circ \exp\right)(\ln t)g(t)\frac{1}{t}\,dt.
$$
Using the change of variable $\tau=\ln t$, we obtain
$$
\int_a^T (J_a^{\sigma} f)(t)g(t)\frac{1}{t}\,dt=\int_{\ln a}^{\ln T} \left(I_{\ln a}^\sigma\,\, f\circ \exp\right)(\tau) (g\circ \exp)(\tau)\,d\tau.
$$
Notice that, since $f\in L^q\left([a,T],\mathbb{R},\frac{1}{t}\,dt\right)$ and $g\in L^r\left([a,T],\mathbb{R},\frac{1}{t}\,dt\right)$, then 
$f\circ \exp \in L^q([\ln a,\ln T],\mathbb{R})$ and $g\circ \exp \in L^r([\ln a,\ln T],\mathbb{R})$. Then, using Lemma \ref{L1}, we obtain 
$$
\int_a^T(J_a^{\sigma} f)(t)g(t)\frac{1}{t}\,dt=\int_{\ln a}^{\ln T} \left(f\circ \exp\right)(\tau) \left(I_{\ln T}^\sigma\,\, g\circ \exp\right)(\tau) \,d\tau.
$$
Using the above change of variable, there holds
$$
\int_a^T(J_a^{\sigma} f)(t)g(t)\frac{1}{t}\,dt=\int_a^T f(t) \left(I_{\ln T}^\sigma\,\, g\circ \exp\right)(\ln t)\frac{1}{t}\,dt.
$$
Thus, by \eqref{in2}, the desired result follows. \hfill $\square$\\

For $\kappa\gg 1$ ($\kappa$ is sufficiently large), let
$$
\mu(t)=\left(\ln \frac{T}{a}\right)^{-\kappa}\left(\ln \frac{T}{t}\right)^{\kappa},\quad  a\leq t\leq T.
$$
\begin{Lemma}\label{L3}
Let $\sigma>0$. Then
\begin{eqnarray}
\label{pp1} (J_T^{\sigma} \mu)(t)&=&\frac{\Gamma(\kappa+1)}{\Gamma(\sigma+\kappa+1)} \left(\ln \frac{T}{a}\right)^{-\kappa}\left(\ln \frac{T}{t}\right)^{\sigma+\kappa},\\
\label{pp2} t(J_T^{\sigma} \mu)'(t)&=&- \frac{\Gamma(\kappa+1)}{\Gamma(\sigma+\kappa)} \left(\ln \frac{T}{a}\right)^{-\kappa}\left(\ln \frac{T}{t}\right)^{\sigma+\kappa-1}.
\end{eqnarray}
\end{Lemma}

\noindent{\it Proof.} We have
\begin{eqnarray*}
(J_T^{\sigma} \mu)(t)&=&\left(\ln \frac{T}{a}\right)^{-\kappa} \frac{1}{\Gamma(\sigma)} \int_t^T  \left(\ln s-\ln t\right)^{\sigma-1} \left(\ln T-\ln s\right)^{\kappa} \frac{1}{s}\,ds\\
&=& \left(\ln \frac{T}{a}\right)^{-\kappa} \frac{1}{\Gamma(\sigma)} \int_t^T  \left(\ln s-\ln t\right)^{\sigma-1} \left((\ln T-\ln t)-(\ln s-\ln t)\right)^{\kappa} \frac{1}{s}\,ds\\
&=& \left(\ln \frac{T}{a}\right)^{-\kappa}  \left(\ln \frac{T}{t}\right)^\kappa \frac{1}{\Gamma(\sigma)}  \int_t^T  \left(\ln s-\ln t\right)^{\sigma-1}\left(1-\frac{\ln s-\ln t}{\ln T-\ln t}\right)^\kappa \frac{1}{s}\,ds.
\end{eqnarray*}
Using the change of variable $\displaystyle\tau=\frac{\ln s-\ln t}{\ln T-\ln t}$, we obtain
\begin{eqnarray*}
(J_T^{\sigma} \mu)(t)&=&	 \left(\ln \frac{T}{a}\right)^{-\kappa}  \left(\ln \frac{T}{t}\right)^{\kappa+\sigma} \frac{1}{\Gamma(\sigma)}  \int_0^1  \tau^{\sigma-1} (1-\tau)^\kappa\,d\tau\\
&=& \left(\ln \frac{T}{a}\right)^{-\kappa}  \left(\ln \frac{T}{t}\right)^{\kappa+\sigma} \frac{1}{\Gamma(\sigma)} B(\sigma,\kappa+1),
\end{eqnarray*}
where $B(\cdot,\cdot)$ is the Beta function. Using the property (see e.g. \cite{KIL})
$$
B(x,y)=\frac{\Gamma(x)\Gamma(y)}{\Gamma(x+y)},\quad x>0,\,y>0,
$$
we get
\begin{eqnarray*}
(J_T^{\sigma} \mu)(t)&=&	 \left(\ln \frac{T}{a}\right)^{-\kappa}  \left(\ln \frac{T}{t}\right)^{\kappa+\sigma} \frac{1}{\Gamma(\sigma)} \frac{\Gamma(\sigma)\Gamma(\kappa+1)}{\Gamma(\sigma+\kappa+1)}\\
&=& \frac{\Gamma(\kappa+1)}{\Gamma(\sigma+\kappa+1)} \left(\ln \frac{T}{a}\right)^{-\kappa}\left(\ln \frac{T}{t}\right)^{\sigma+\kappa},
\end{eqnarray*}
which proves \eqref{pp1}.   Differentiating \eqref{pp1} and using the property (see e.g. \cite{KIL})
$$
x\Gamma(x)=\Gamma(x+1),\quad x>0,
$$
\eqref{pp2} follows. \hfill $\square$\\

Let  $f\in AC([a,T],\mathbb{R})$. The left/forward Hadamard-Caputo fractional derivative of order $\alpha\in (0,1)$ of $f$, is defined by (see Agrawal \cite{AG})
\begin{eqnarray*}
(\mathcal{D}_{a}^\alpha f)(t)&=&J_a^{1-\alpha} \left(t f'\right)(t)\\
&=& \frac{1}{\Gamma(1-\alpha)} \int_a^t  \left(\ln \frac{t}{s}\right)^{-\alpha} f'(s)\,ds,
\end{eqnarray*}
for almost everywhere $t\in [a,T]$. 

Let $u: [a,T]\times \mathbb{R}^N\to \mathbb{C}$ be a given complex-valued function.  For a fixed $x\in \mathbb{R}^N$, we denote by $u(\cdot,x): [a,T]\to \mathbb{C}$ the function defined by
$$
u(\cdot,x)(t)=u(t,x),\quad t\in [a,T].
$$
The left-sided  and right-sided Hadamard fractional integrals of  order $\sigma>0$ of  $u$ with respect to the time variable $t$, are defined respectively by 
\begin{eqnarray*}
(J_{a|t}^\sigma u)(t,x)&=&(J_a^\sigma\,\, u(\cdot,x))(t)\\
&=&(J_a^\sigma\,\, \Re u(\cdot,x))(t)+i (J_a^\sigma\,\, \Im u(\cdot,x))(t) \\
&=& \frac{1}{\Gamma(\sigma)} \int_a^t  \left(\ln \frac{t}{s}\right)^{\sigma-1} \Re u(s,x) \frac{1}{s}\,ds+\frac{i}{\Gamma(\sigma)} \int_a^t  \left(\ln \frac{t}{s}\right)^{\sigma-1} \Im u(s,x) \frac{1}{s}\,ds
\end{eqnarray*}
and
\begin{eqnarray*}
(J_{T|t}^\sigma u)(t,x)&=&(J_T^\sigma\,\, u(\cdot,x))(t)\\
&=&(J_T^\sigma\,\, \Re u(\cdot,x))(t)+i (J_T^\sigma\,\, \Im u(\cdot,x))(t)\\
&=& \frac{1}{\Gamma(\sigma)} \int_t^T  \left(\ln \frac{s}{t}\right)^{\sigma-1} \Re u(s,x) \frac{1}{s}\,ds+\frac{i}{\Gamma(\sigma)} \int_t^T  \left(\ln \frac{s}{t}\right)^{\sigma-1} \Im u(s,x) \frac{1}{s}\,ds.
\end{eqnarray*}
The left/forward Hadamard-Caputo fractional derivative of order $\alpha\in (0,1)$ of $u$ with respect to the time variable $t$, is defined by
\begin{eqnarray*}
(\mathcal{D}_{a|t}^\alpha u)(t,x)&=&J_{a|t}^{1-\alpha} \left(t \frac{\partial u}{\partial t}(\cdot,x)\right)(t)	\\
&=& \frac{1}{\Gamma(1-\alpha)} \int_a^t  \left(\ln \frac{t}{s}\right)^{-\alpha} \frac{\partial \Re u}{\partial s}(s,x)\,ds+\frac{i}{\Gamma(1-\alpha)} \int_a^t  \left(\ln \frac{t}{s}\right)^{-\alpha} \frac{\partial \Im u}{\partial s}(s,x)\,ds. 
\end{eqnarray*}

\section{Main results}\label{sec3}

Before stating our main results, let us give the definition of global weak solutions to \eqref{P}. Observe that \eqref{P} is equivalent to the system 
\begin{eqnarray}\label{Sys-P}
\left\{\begin{array}{llll}
r_\alpha \mathcal{D}_{a|t}^\alpha u_1-s_\alpha \mathcal{D}_{a|t}^\alpha u_2+\Delta u_1 = \lambda_1 \displaystyle \left(\ln \frac{t}{a}\right)^\gamma |u|^p,\quad t>a,\,\, x\in \mathbb{R}^N,\\ \\
s_\alpha \mathcal{D}_{a|t}^\alpha u_1+r_\alpha  \mathcal{D}_{a|t}^\alpha u_2+\Delta u_2 = \lambda_2 \displaystyle \left(\ln \frac{t}{a}\right)^\gamma |u|^p,\quad t>a,\,\, x\in \mathbb{R}^N,\\ \\
(u_1(a,x),u_2(a,x))= (f_1(x),f_2(x)),\quad x\in \mathbb{R}^N,
\end{array}
\right.	
\end{eqnarray}
where 
$$
(r_\alpha,s_\alpha)=\left(\cos\left(\frac{\pi\alpha}{2}\right),\sin\left(\frac{\pi\alpha}{2}\right)\right),\,\, (u_1,u_2)=(\Re u,\Im u),\,\,(f_1,f_2)=(\Re f,\Im f),\,\, (\lambda_1,\lambda_2)=(\Re \lambda,\Im \lambda).
$$

For $T>0$, let
$$
Q_T=[a,T]\times \mathbb{R}^N
$$
and $\Phi_T$ be the set of functions $\varphi$ satisfying:
$$
\varphi\in C_{t,x}^{1,2}	(Q_T,\mathbb{R}),\quad \mbox{supp}_x\varphi\subset\subset \mathbb{R}^N.
$$
Multiplying the first two equations in \eqref{Sys-P} by $\varphi\in \Phi_T$, using the integration by parts rule provided by Lemma \ref{L2},  the initial conditions in \eqref{Sys-P}, and taking in consideration Remark \ref{RK-PR}, we  define global weak solutions to \eqref{P} as follows. 

\begin{Definition}\label{WS-Sys-P}
We say that $u$ is a global weak solution to \eqref{P}, if 
\begin{itemize}
\item[{\rm{(i)}}] 
$$
u\in L^1_{loc}([a,\infty)\times \mathbb{R}^N,\mathbb{C}),\quad \left(\ln \frac{t}{a}\right)^\gamma |u|^p\in L^1_{loc}([a,\infty)\times \mathbb{R}^N,\mathbb{R}),
$$
\item[{\rm{(ii)}}] for all $T>0$ and $\varphi\in \Phi_T$, there holds
\begin{equation}\label{ws1}
\begin{aligned}
&\lambda_1\int_{Q_T}	\left(\ln \frac{t}{a}\right)^\gamma |u|^p\varphi\,dx\,dt+\int_{\mathbb{R}^N}\left(r_\alpha f_1(x)-s_\alpha f_2(x)\right)(J_{T|t}^{1-\alpha}\,\,t\varphi)(a,x)\,dx\\
&=\int_{Q_T} u_1\Delta \varphi\,dx\,dt-\int_{Q_T} \left(r_\alpha u_1-s_\alpha u_2\right)\frac{\partial J_{T|t}^{1-\alpha}\,\,t \varphi}{\partial t}\,dx\,dt
\end{aligned}
\end{equation}
and
\begin{equation}\label{ws2}
\begin{aligned}
&\lambda_2\int_{Q_T}	\left(\ln \frac{t}{a}\right)^\gamma |u|^p\varphi\,dx\,dt+\int_{\mathbb{R}^N}\left(s_\alpha f_1(x)+r_\alpha f_2(x)\right)(J_{T|t}^{1-\alpha}\,\,t\varphi)(a,x)\,dx\\
&=\int_{Q_T} u_2\Delta \varphi\,dx\,dt-\int_{Q_T} \left(s_\alpha u_1+r_\alpha u_2\right)\frac{\partial J_{T|t}^{1-\alpha}\,\,t \varphi}{\partial t}\,dx\,dt.
\end{aligned}
\end{equation}
\end{itemize}
\end{Definition}

Our  main results are the following.

\begin{Theorem}\label{T1}
Let $\gamma>-\alpha$, $\gamma(N\alpha-2)<2\alpha$, and 
\begin{equation}\label{blow-up1}
\max\{1,1+\gamma\}<p<1+\frac{2(\alpha+\gamma)}{N\alpha}.	
\end{equation}
If the initial value $f\in L^1(\mathbb{R}^N,\mathbb{C})$  satisfies
\begin{equation}\label{assump-f}
\lambda_1 	\int_{\mathbb{R}^N}\left(r_\alpha f_1(x)-s_\alpha f_2(x)\right)\,dx>0\quad\mbox{or}\quad  \lambda_2 	\int_{\mathbb{R}^N}\left(s_\alpha f_1(x)+r_\alpha f_2(x)\right)\,dx>0,
\end{equation} 
 then \eqref{P} admits no global weak solution. 	
\end{Theorem}

\begin{Remark} 
Notice that under the conditions $\gamma>-\alpha$ and  $\gamma(N\alpha-2)<2\alpha$, the set of $p$ satisfying \eqref{blow-up1} is nonempty.
\end{Remark}

\begin{Theorem}\label{T2}
Let $\gamma>0$ and 
\begin{equation}\label{blow-upT2}
1+\gamma<p<1+\frac{\gamma}{\alpha}.	
\end{equation}	
If the initial value $f\in L^1(\mathbb{R}^N,\mathbb{C})$  satisfies \eqref{assump-f}, then \eqref{P} admits no global weak solution. 	
\end{Theorem}

In the case $\gamma>0$,  we deduce from Theorems \ref{T1}and \ref{T2} the following result. 

\begin{Corollary}\label{CR1}
Let $\gamma>0$, $\gamma(N\alpha-2)<2\alpha$, and 
$$
1+\gamma<p< \max\left\{1+\frac{2(\alpha+\gamma)}{N\alpha},1+\frac{\gamma}{\alpha}\right\}.	
$$
If the initial value $f\in L^1(\mathbb{R}^N,\mathbb{C})$  satisfies \eqref{assump-f}, then \eqref{P} admits no global weak solution. 		
\end{Corollary}

\begin{Remark}\label{RK-Cor}
{\rm{(i)}} Notice that for $\gamma>0$, we have
$$
\max\left\{1+\frac{2(\alpha+\gamma)}{N\alpha},1+\frac{\gamma}{\alpha}\right\}=\left\{\begin{array}{llll}
\displaystyle 1+\frac{2(\alpha+\gamma)}{N\alpha} &\mbox{if}& (N-2)\gamma<2\alpha,\\ \\
\displaystyle 1+\frac{\gamma}{\alpha} &\mbox{if}& \gamma(N\alpha-2)<2\alpha \leq (N-2)\gamma.	
\end{array}
\right.
$$
{\rm{(ii)}} Observe that, if $N\in\{1,2\}$ and $\gamma>0$, then   $(N-2)\gamma<2\alpha$. Hence, by {\rm{(i)}}, we deduce that 
$$
\max\left\{1+\frac{2(\alpha+\gamma)}{N\alpha},1+\frac{\gamma}{\alpha}\right\}=1+\frac{2(\alpha+\gamma)}{N\alpha}.
$$ 
\end{Remark}

We provide below some examples to illustrate our obtained results.

\begin{Example}
Consider the initial value problem for the nonlinear time-fractional  Schr\"odinger equation 
\begin{eqnarray}\label{P-ex1}
\left\{\begin{array}{llll}
\sqrt{i} \mathcal{D}_{a|t}^{1/2} u+\Delta u	=\displaystyle \left(\ln \frac{t}{a}\right)^{-1/4}  |u|^p,\quad t>a,\,\, x\in \mathbb{R}^N,\\ \\
u(a,x)=\displaystyle\frac{1}{|x|^{N-1}(1+|x|^2)},\quad x\in \mathbb{R}^N,
\end{array}
\right.
\end{eqnarray}
where $a>0$ and $N\geq 5$. Then \eqref{P-ex1} is a special case of \eqref{P} with 
$$
\alpha=\frac{1}{2},\quad \gamma=-\frac{1}{4},\quad \lambda=\lambda_1=1,\quad f(x)=f_1(x)=\frac{1}{|x|^{N-1}(1+|x|^2)}.
$$
Observe that $f\in L^1(\mathbb{R}^N,\mathbb{R})$ and 
$$
\lambda_1 	\int_{\mathbb{R}^N}\left(r_\alpha f_1(x)-s_\alpha f_2(x)\right)\,dx=\frac{\sqrt{2}}{2}\int_{\mathbb{R}^N} \frac{1}{|x|^{N-1}(1+|x|^2)}\,dx>0.
$$
Moreover, we have 
$$
\gamma>-\alpha,\quad \gamma(N\alpha-2)=-\frac{1}{8}(N-4)<0<2\alpha,\quad \max\{1,1+\gamma\}=1,\quad 1+\frac{2(\alpha+\gamma)}{N\alpha} =1+\frac{1}{N}.
$$
Hence, by Theorem \ref{T1}, we deduce that for all
$$
1<p<1+\frac{1}{N},
$$
\eqref{P-ex1} admits no global weak solution.
\end{Example}

\begin{Example}
Consider the initial value problem for the nonlinear time-fractional  Schr\"odinger equation 
\begin{eqnarray}\label{P-ex2}
\left\{\begin{array}{llll}
\sqrt{i} \mathcal{D}_{a|t}^{1/2} u+\Delta u	=-\displaystyle \left(\ln \frac{t}{a}\right)^{1/N}  |u|^p,\quad t>a,\,\, x\in \mathbb{R}^N,\\ \\
u(a,x)=i|x|^{2-N}\exp(-|x|^2),\quad x\in \mathbb{R}^N,
\end{array}
\right.
\end{eqnarray}	
where $a>0$ and $N\geq 1$. Then \eqref{P-ex2} is a special case of \eqref{P} with 
$$
\alpha=\frac{1}{2},\quad \gamma=\frac{1}{N},\quad \lambda=\lambda_1=-1,\quad f(x)=if_2(x)=i|x|^{2-N}\exp(-|x|^2).
$$
Observe that $f\in L^1(\mathbb{R}^N,\mathbb{C})$ and  
$$
\lambda_1 	\int_{\mathbb{R}^N}\left(r_\alpha f_1(x)-s_\alpha f_2(x)\right)\,dx=\frac{\sqrt{2}}{2}\int_{\mathbb{R}^N} |x|^{2-N}\exp(-|x|^2)\,dx>0.
$$
Moreover, we have
$$
\gamma>0,\quad (N-2)\gamma=1-\frac{2}{N}<1=2\alpha,\quad \max\{1,1+\gamma\}=1+\frac{1}{N},\quad 1+\frac{2(\alpha+\gamma)}{N\alpha}=1+\frac{2}{N}+\frac{4}{N^2}.
$$
Hence by Corollary \ref{CR1}, and taking in consideration Remark \ref{RK-Cor} {\rm{(i)}}, we deduce that for all
$$
1+\frac{1}{N}<p< 1+\frac{2}{N}+\frac{4}{N^2}, 
$$
\eqref{P-ex2} admits no global weak solution.
\end{Example}

\begin{Example}
Consider the initial value problem for the nonlinear time-fractional  Schr\"odinger equation 
\begin{eqnarray}\label{P-ex3}
\left\{\begin{array}{llll}
\sqrt{i} \mathcal{D}_{a|t}^{1/2} u+\Delta u	=\displaystyle \left(\ln \frac{t}{a}\right)^{1/(N-2)}  |u|^p,\quad t>a,\,\, x\in \mathbb{R}^N,\\ \\
u(a,x)=\exp(-|x|),\quad x\in \mathbb{R}^N,
\end{array}
\right.
\end{eqnarray}	
where $a>0$ and $N\geq 3$. Then \eqref{P-ex3} is a special case of \eqref{P} with 
$$
\alpha=\frac{1}{2},\quad \gamma=\frac{1}{N-2},\quad \lambda=\lambda_1=1,\quad f(x)=f_1(x)=\exp(-|x|).
$$
Observe that $f\in L^1(\mathbb{R}^N,\mathbb{R})$ and  
$$
\lambda_1 	\int_{\mathbb{R}^N}\left(r_\alpha f_1(x)-s_\alpha f_2(x)\right)\,dx=\frac{\sqrt{2}}{2}\int_{\mathbb{R}^N} \exp(-|x|)\,dx>0.
$$
Moreover, we have
$$
\gamma>0,\quad \gamma(N\alpha-2)=\frac{N-4}{2(N-2)}<2\alpha=(N-2)\gamma,\quad 1+\frac{\gamma}{\alpha}=1+\frac{2}{N-2}.
$$
Hence by Corollary \ref{CR1}, and taking in consideration Remark \ref{RK-Cor} {\rm{(i)}}, we deduce that for all
$$
1+\frac{1}{N-2}<p< 1+\frac{2}{N-2}, 
$$
\eqref{P-ex3} admits no global weak solution.
\end{Example}

\section{Proof of the main results}\label{sec4}

In the sequel,  we use $C$ to denote a positive constant which may vary from line to line, but its value  is not essential to the analysis of the problem. The proof of our main results is based on the test function method developed by Mitidieri and Pohozaev \cite{MI}, and a judicious choice of the test function. 

\noindent{\it Proof of Theorem \ref{T1}.}  We argue by contradiction. Namely, we suppose that $u$ is a global weak solution to \eqref{P}.  Let 
\begin{equation}\label{lambda1I}
\lambda_1 	\int_{\mathbb{R}^N}\left(r_\alpha f_1(x)-s_\alpha f_2(x)\right)\,dx>0.
\end{equation}
Then, by \eqref{ws1}, for all $T>0$ and $\varphi\in \Phi_T$, we have (after a multiplication by $\lambda_1$)
$$
\begin{aligned}
&\lambda_1^{2}\int_{Q_T}	\left(\ln \frac{t}{a}\right)^\gamma |u|^p\varphi\,dx\,dt+\lambda_1\int_{\mathbb{R}^N}\left(r_\alpha f_1(x)-s_\alpha f_2(x)\right)(J_{T|t}^{1-\alpha}\,\,t\varphi)(a,x)\,dx\\
&=\lambda_1 \int_{Q_T} u_1\Delta \varphi\,dx\,dt-\lambda_1 \int_{Q_T} \left(r_\alpha u_1-s_\alpha u_2\right)\frac{\partial J_{T|t}^{1-\alpha}\,\,t \varphi}{\partial t}\,dx\,dt,
\end{aligned}
$$
which yields
\begin{equation}\label{gest1}
\begin{aligned}
&\lambda_1^2\int_{Q_T}	\left(\ln \frac{t}{a}\right)^\gamma |u|^p\varphi\,dx\,dt+\lambda_1\int_{\mathbb{R}^N}\left(r_\alpha f_1(x)-s_\alpha f_2(x)\right)(J_{T|t}^{1-\alpha}\,\,t\varphi)(a,x)\,dx\\
&\leq|\lambda_1|\int_{Q_T}|u| |\Delta \varphi|\,dx\,dt+2|\lambda_1|\int_{Q_T} |u|\left|\frac{\partial J_{T|t}^{1-\alpha}\,\,t \varphi}{\partial t}\right|\,dx\,dt.
\end{aligned}
\end{equation}
On the other hand, by $\varepsilon$-Young inequality with $\varepsilon=\frac{|\lambda_1|}{2}>0$, we have
\begin{equation}\label{est1}
\int_{Q_T}|u| |\Delta \varphi|\,dx\,dt \leq \frac{|\lambda_1|}{2}\int_{Q_T}	\left(\ln \frac{t}{a}\right)^\gamma |u|^p\varphi\,dx\,dt+C  \int_{Q_T}	\left(\ln \frac{t}{a}\right)^{-\frac{\gamma}{p-1}}\varphi^{-\frac{1}{p-1}} |\Delta \varphi|^{\frac{p}{p-1}}\,dx\,dt.
\end{equation}
Similarly, using $\varepsilon$-Young inequality with $\varepsilon=\frac{|\lambda_1|}{4}$, we get
\begin{equation}\label{est2}
\begin{aligned}
&\int_{Q_T} |u|\left|\frac{\partial J_{T|t}^{1-\alpha}\,\,t \varphi}{\partial t}\right|\,dx\,dt\\
&\leq \frac{|\lambda_1|}{4}\int_{Q_T}	\left(\ln \frac{t}{a}\right)^\gamma |u|^p\varphi\,dx\,dt+C  \int_{Q_T}	\left(\ln \frac{t}{a}\right)^{-\frac{\gamma}{p-1}}\varphi^{-\frac{1}{p-1}} \left|\frac{\partial J_{T|t}^{1-\alpha}\,\,t \varphi}{\partial t}\right|^{\frac{p}{p-1}}\,dx\,dt.	
\end{aligned}
\end{equation}
Hence, combining \eqref{gest1}, \eqref{est1}, and \eqref{est2}, we obtain
\begin{equation}\label{estJ1J2}
\lambda_1\int_{\mathbb{R}^N}\left(r_\alpha f_1(x)-s_\alpha f_2(x)\right)(J_{T|t}^{1-\alpha}\,\,t\varphi)(a,x)\,dx\leq C(K_1(\varphi)+K_2(\varphi)),	
\end{equation}
where
$$
K_1(\varphi)=\int_{Q_T}	\left(\ln \frac{t}{a}\right)^{-\frac{\gamma}{p-1}}\varphi^{-\frac{1}{p-1}} \left|\frac{\partial J_{T|t}^{1-\alpha}\,\,t \varphi}{\partial t}\right|^{\frac{p}{p-1}}\,dx\,dt
$$
and
$$
K_2(\varphi)=	\int_{Q_T}	\left(\ln \frac{t}{a}\right)^{-\frac{\gamma}{p-1}}\varphi^{-\frac{1}{p-1}} |\Delta \varphi|^{\frac{p}{p-1}}\,dx\,dt.
$$
Now, consider a family of cut-off functions $\{\xi_R\}_{R\geq 1}\subset C_c^\infty(\mathbb{R}^N,\mathbb{R})$ ($\xi_R\in C^\infty(\mathbb{R}^N,\mathbb{R})$ and $\mbox{supp }\xi_R\subset\subset \mathbb{R}^N$) satisfying the following properties:
\begin{itemize}
\item[(a)] 	$0\leq \xi_R\leq 1$, ${\xi_R}{|_{B_R}}\equiv 1$,
\item[(b)] $\mbox{\rm{supp}}(\xi_R)\subset B_{2R}$,
\item[(c)] $|\nabla \xi_R|\leq \frac{C}{R}$,
\item[(d)] $|\Delta \xi_R|\leq \frac{C}{R^2}$,
\end{itemize}
where for $\rho>0$, 
$$
B_\rho=\{x\in \mathbb{R}^N:\, |x|<\rho\}.
$$
For $\kappa\gg 1$ and $\ell\gg1 $, let us introduce the test function 
\begin{equation}\label{testf-choice}
\varphi(t,x)=\eta(t) \xi_R^\ell(x),\quad a\leq t\leq T, 	
\end{equation}
where 
$$
\eta(t)=\frac{1}{t} \left(\ln \frac{T}{a}\right)^{-\kappa}\left(\ln \frac{T}{t}\right)^{\kappa}=\frac{1}{t}\mu(t).
$$
It can be easily seen that for all $T>0$, the function $\varphi$ defined by \eqref{testf-choice}
belongs to $\Phi_T$, and thus, it satisfies the estimate \eqref{estJ1J2}.  

Let us estimate the terms $K_j(\varphi)$, $j=1,2$, where $\varphi$ is defined by \eqref{testf-choice}. The term $K_1(\varphi)$ can be written as 
\begin{equation}\label{K1term}
K_1(\varphi)=K_{11}(\varphi)K_{12}(\varphi),
\end{equation}
where
$$
K_{11}(\varphi)=\int_a^T \left(\ln \frac{t}{a}\right)^{-\frac{\gamma}{p-1}}t^{\frac{1}{p-1}}\mu^{-\frac{1}{p-1}}(t)\left|(J_{T}^{1-\alpha}\,\,\mu)'(t)\right|^{\frac{p}{p-1}}\,dt
$$
and
$$
K_{12}(\varphi)=\int_{\mathbb{R}^N} \xi_R^\ell(x)\,dx.
$$
On the other hand, using \eqref{pp2}, for $a<t<T$, we get 
$$
\mu^{-\frac{1}{p-1}}(t)\left|(J_{T}^{1-\alpha}\,\,\mu)'(t)\right|^{\frac{p}{p-1}}=\left[\frac{\Gamma(\kappa+1)}{\Gamma(1-\alpha+\kappa)}\right]^{\frac{p}{p-1}}\left(\ln \frac{T}{a}\right)^{-\kappa}\left(\ln \frac{T}{t}\right)^{\kappa-\frac{\alpha p}{p-1}}t^{-\frac{p}{p-1}},
$$
which yields
\begin{eqnarray*}
K_{11}(\varphi)&=&\left[\frac{\Gamma(\kappa+1)}{\Gamma(1-\alpha+\kappa)}\right]^{\frac{p}{p-1}}\left(\ln \frac{T}{a}\right)^{-\kappa}\int_a^T\left(\ln \frac{t}{a}\right)^{-\frac{\gamma}{p-1}} \left(\ln \frac{T}{t}\right)^{\kappa-\frac{\alpha p}{p-1}}\frac{1}{t}\,dt\\
&\leq & \left[\frac{\Gamma(\kappa+1)}{\Gamma(1-\alpha+\kappa)}\right]^{\frac{p}{p-1}}\left(\ln \frac{T}{a}\right)^{-\frac{\alpha p}{p-1}}\int_a^T\left(\ln \frac{t}{a}\right)^{-\frac{\gamma}{p-1}} \frac{1}{t}\,dt.
\end{eqnarray*}
Then, since $\gamma<p-1$ by \eqref{blow-up1}, there holds
\begin{equation}\label{K11est}
K_{11}(\varphi)\leq \left[\frac{\Gamma(\kappa+1)}{\Gamma(1-\alpha+\kappa)}\right]^{\frac{p}{p-1}}\left(1-\frac{\gamma}{p-1}\right)^{-1}	\left(\ln \frac{T}{a}\right)^{1-\frac{\gamma+\alpha p}{p-1}}.
\end{equation}
Next, by the properties (a) and (b) of $\xi_R$, we get
\begin{eqnarray}\label{K12est}
\nonumber K_{12}(\varphi)&=&\int_{B_{2R}}\xi_R^\ell(x)\,dx\\
\nonumber &\leq & \mbox{Vol}(B_{2R})\\
&=& CR^N.	
\end{eqnarray}
Hence, it follows from \eqref{K1term}, \eqref{K11est}, and \eqref{K12est} that
\begin{equation}\label{estK1phi}
K_1(\varphi)\leq C 	\left(\ln \frac{T}{a}\right)^{1-\frac{\gamma+\alpha p}{p-1}}R^N.	
\end{equation}
The term $K_2(\varphi)$ can be written as 
\begin{equation}\label{K2term}
K_2(\varphi)=K_{21}(\varphi)K_{22}(\varphi),
\end{equation}
where
$$
K_{21}(\varphi)=\int_a^T \left(\ln \frac{t}{a}\right)^{-\frac{\gamma}{p-1}}\eta(t)\,dt
$$
and
$$
K_{22}(\varphi)=\int_{\mathbb{R}^N}\xi_R^{-\frac{\ell}{p-1}}(x)\left|\Delta[\xi_R^\ell(x)]\right|^{\frac{p}{p-1}}\,dx.
$$
On the other hand, we have
\begin{eqnarray}\label{K21est}
\nonumber K_{21}(\varphi)	&=& \left(\ln \frac{T}{a}\right)^{-\kappa}\int_a^T \left(\ln \frac{t}{a}\right)^{-\frac{\gamma}{p-1}}\frac{1}{t} \left(\ln \frac{T}{t}\right)^{\kappa}\\
\nonumber &\leq & \int_a^T \left(\ln \frac{t}{a}\right)^{-\frac{\gamma}{p-1}}\frac{1}{t} \,dt\\
&=& \left(1-\frac{\gamma}{p-1}\right)^{-1}\left(\ln \frac{T}{a}\right)^{1-\frac{\gamma}{p-1}}
\end{eqnarray}
Next, using the property
$$
\Delta(\xi_R^\ell)=\ell\xi_R^{\ell-2}\left((\ell-1)|\nabla \xi_R|^2+\xi_R\Delta \xi_R\right)
$$
as well as the properties  (a)--(d) of $\xi_R$, we obtain
$$
\left|\Delta[\xi_R^\ell(x)]\right|\leq C R^{-2}\xi_R^{\ell-2},\quad R<|x|<2R,
$$
which yields
\begin{eqnarray}\label{est-K22}
\nonumber K_{22}(\varphi)	&\leq & C R^{-\frac{2p}{p-1}}\int_{R<|x|<2R} \xi_R^{\ell-\frac{2p}{p-1}}(x)\,dx\\
\nonumber &\leq & C  R^{-\frac{2p}{p-1}}\mbox{Vol}(B_{2R})\\
&=& C R^{N-\frac{2p}{p-1}}.
\end{eqnarray}
Hence, by \eqref{K2term}, \eqref{K21est}, and \eqref{est-K22}, we deduce that 
\begin{equation}\label{est-K2}
K_2(\varphi)\leq C 	\left(\ln \frac{T}{a}\right)^{1-\frac{\gamma}{p-1}} R^{N-\frac{2p}{p-1}}.
\end{equation}
Now, consider the term from the left-hand side of \eqref{estJ1J2}.  By \eqref{testf-choice} and using \eqref{pp1}, we have
\begin{eqnarray*}
(J_{T|t}^{1-\alpha}\,\,t\varphi)(a,x)&=&(J_{T|t}^{1-\alpha}\,\,t\eta \xi_R^\ell)(a,x)\\
&=& \frac{\Gamma(\kappa+1)}{\Gamma(\kappa+2-\alpha)}\left(\ln \frac{T}{a}\right)^{1-\alpha}\xi_R^\ell(x),
\end{eqnarray*}
which yields
\begin{equation}\label{right-term}
\begin{aligned}
&\lambda_1\int_{\mathbb{R}^N}\left(r_\alpha f_1(x)-s_\alpha f_2(x)\right)(J_{T|t}^{1-\alpha}\,\,t\varphi)(a,x)\,dx\\
&=\frac{\Gamma(\kappa+1)}{\Gamma(\kappa+2-\alpha)}	\left(\ln \frac{T}{a}\right)^{1-\alpha}\lambda_1 \int_{\mathbb{R}^N}\left(r_\alpha f_1(x)-s_\alpha f_2(x)\right)\xi_R^\ell(x)\,dx.
\end{aligned}
\end{equation}
Then, it follows from \eqref{estJ1J2}, \eqref{estK1phi}, \eqref{est-K2}, and \eqref{right-term} that 
\begin{equation}\label{finest-case1}
\lambda_1 \int_{\mathbb{R}^N}\left(r_\alpha f_1(x)-s_\alpha f_2(x)\right)\xi_R^\ell(x)\,dx\leq C  \left(\left(\ln T\right)^{\alpha-\frac{\gamma+\alpha p}{p-1}}R^N+ \left(\ln T\right)^{\alpha-\frac{\gamma}{p-1}} R^{N-\frac{2p}{p-1}}\right),\quad T\gg 1.
\end{equation}
Next, for sufficiently large $R$, taking $T=\exp R^\theta$, where $\theta>0$ is a constant that will be chosen later, \eqref{finest-case1} reduces to 
$$
\lambda_1 \int_{\mathbb{R}^N}\left(r_\alpha f_1(x)-s_\alpha f_2(x)\right)\xi_R^\ell(x)\,dx\leq C \left(R^{N+\theta\left(\alpha-\frac{\gamma+\alpha p}{p-1}\right)}+R^{N-\frac{2p}{p-1}+\theta\left(\alpha-\frac{\gamma}{p-1}\right)}\right).
$$
On the other hand, observe that for $\theta=\frac{2}{\alpha}$, we have
$$
N+\theta\left(\alpha-\frac{\gamma+\alpha p}{p-1}\right)=N-\frac{2p}{p-1}+\theta\left(\alpha-\frac{\gamma}{p-1}\right)=\frac{N\alpha(p-1)-2(\alpha+\gamma)}{\alpha(p-1)}.
$$
Hence, for this value of $\theta$, we have
\begin{equation}\label{eqlim}
\lambda_1 \int_{\mathbb{R}^N}\left(r_\alpha f_1(x)-s_\alpha f_2(x)\right)\xi_R^\ell(x)\,dx\leq C R^{\frac{N\alpha(p-1)-2(\alpha+\gamma)}{\alpha (p-1)}}. 	
\end{equation}
Notice that by \eqref{blow-up1}, we have
$$
\frac{N\alpha(p-1)-2(\alpha+\gamma)}{\alpha(p-1)}<0.
$$
Hence, since $f\in L^1(\mathbb{R}^N,\mathbb{C})$, using the property (a) of the cut-off function $\xi_R$, the dominated convergence theorem, and passing to the limit as $R\to \infty$ in \eqref{eqlim}, we obtain
$$
\lambda_1 \int_{\mathbb{R}^N}\left(r_\alpha f_1(x)-s_\alpha f_2(x)\right)\,dx\leq 0,
$$ 
which contradicts \eqref{lambda1I}.

Suppose now that 
\begin{equation}\label{lambda2pos}
\lambda_2 	\int_{\mathbb{R}^N}\left(s_\alpha f_1(x)+r_\alpha f_2(x)\right)\,dx>0.
\end{equation}
Then, by \eqref{ws2},  we have (after a multiplication by $\lambda_2$)
$$
\begin{aligned}
&\lambda_2^2\int_{Q_T}	\left(\ln \frac{t}{a}\right)^\gamma |u|^p\varphi\,dx\,dt+\lambda_2\int_{\mathbb{R}^N}\left(s_\alpha f_1(x)+r_\alpha f_2(x)\right)(J_{T|t}^{1-\alpha}\,\,t\varphi)(a,x)\,dx\\
&=\lambda_2\int_{Q_T} u_2\Delta \varphi\,dx\,dt-\lambda_2\int_{Q_T} \left(s_\alpha u_1+r_\alpha u_2\right)\frac{\partial J_{T|t}^{1-\alpha}\,\,t \varphi}{\partial t}\,dx\,dt,
\end{aligned}
$$
which yields
$$
\begin{aligned}
&\lambda_2^2\int_{Q_T}	\left(\ln \frac{t}{a}\right)^\gamma |u|^p\varphi\,dx\,dt+\lambda_2\int_{\mathbb{R}^N}\left(s_\alpha f_1(x)+r_\alpha f_2(x)\right)(J_{T|t}^{1-\alpha}\,\,t\varphi)(a,x)\,dx\\
&\leq|\lambda_2|\int_{Q_T}|u| |\Delta \varphi|\,dx\,dt+2|\lambda_2|\int_{Q_T} |u|\left|\frac{\partial J_{T|t}^{1-\alpha}\,\,t \varphi}{\partial t}\right|\,dx\,dt.
\end{aligned}
$$
Then, repeating the same procedures as in the previous case, we arrive at
$$
\lambda_2 \int_{\mathbb{R}^N}\left(s_\alpha f_1(x)+r_\alpha f_2(x)\right)\xi_R^\ell(x)\,dx\leq C  \left(\left(\ln T\right)^{\alpha-\frac{\gamma+\alpha p}{p-1}}R^N+ \left(\ln T\right)^{\alpha-\frac{\gamma}{p-1}} R^{N-\frac{2p}{p-1}}\right).
$$
Following exactly   the same steps as above, we reach a contradiction with \eqref{lambda2pos}. The proof of Theorem \ref{T1} is completed. \hfill $\square$\\

\noindent{\it Proof of Theorem \ref{T2}.}   Suppose that $u$ is a global weak solution to \eqref{P}.   We only consider the case \eqref{lambda1I}, since the case  
\eqref{lambda2pos} can be treated in the same way. From the proof of Theorem \ref{T1}, \eqref{finest-case1} holds for $T\gg 1$. Observe that for $\gamma>0$ and under condition \eqref{blow-upT2}, we have
$$
\alpha-\frac{\gamma+\alpha p}{p-1}=-\frac{\alpha+\gamma}{p-1}<0,\quad \alpha-\frac{\gamma}{p-1}<0.
$$
Hence, fixing $R$, and passing to the limit as $T\to \infty$ in \eqref{finest-case1}, we obtain 
$$
\lambda_1 \int_{\mathbb{R}^N}\left(r_\alpha f_1(x)-s_\alpha f_2(x)\right)\xi_R^\ell(x)\,dx\leq 0.
$$
Next, passing to the limit as $R\to \infty$ in the above inequality, we reach a contradiction with \eqref{lambda1I}. The proof of Theorem \ref{T2} is completed. \hfill $\square$

\subsection*{Acknowledgments}
The third author wish to thank  GNAMPA 2021 and the RUDN University Strategic Academic Leadership Program.\\
The fourth author is supported by Researchers Supporting Project number (RSP--2021/4), King Saud University, Riyadh, Saudi Arabia.


\begin{thebibliography}{99} 

\bibitem{AB}
M.I. Abbas, M.A. Ragusa,  On the hybrid fractional differential equations with fractional proportional derivatives of a function with respect to a certain function, Symmetry. 13(2) (2021) 264.

\bibitem{AG}
O.P. Agrawal, Generalized multiparameters fractional variational calculus, Int J Differ Equ. 2012 (2012) 521750.

\bibitem{AR}
O.A. Arqub, Application of residual power series method for the solution of time-fractional schr\"odinger equations in one-dimensional space, Fundam Inform. 166 (2) (2019) 87--110.

\bibitem{BAG}
R.L. Bagley, P.L. Torvik, On the fractional calculus models of viscoelastic behaviour, 
J. Rheol. 30 (1986) 133--155.

\bibitem{DO}
J. Dong, M. Xu, Solutions to the space fractional Schr\"odinger equation using momentum representation method, J Math Phys. 48 (2007) 072105.

\bibitem{Fujita}
H. Fujita, On the blowing up of solutions to the Cauchy problem for $u_t=\Delta u+u^{1+\alpha}$, J. Fac. Sci. Univ. Tokyo, Sect. 1A, Math. 13 (1966) 199--124.

\bibitem{GU}
X. Guo, M. Xu, Some physical applications of fractional Schr\"odinger equation,
J Math Phys. 47 (2006)  82104.

\bibitem{HI}
R. Hilfer, Applications of fractional calculus in physics, World Scientific, New Jersey, 2001.

\bibitem{IKI}
M. Ikeda, T. Inui, Small data blow-up of $L^2$ or $H^1$-solution for the semilinear Schr\"dinger equation without gauge invariance, J. Evol. Equ. 15 (2015) 1--11.


\bibitem{IK}
M. Ikeda, Y. Wakasugi, Small data blow-up of $L^2$-solution for the nonlinear Schr\"dinger equation without gauge invariance, Diff. Int. Equ. 26 (2013) 1275--1285.


\bibitem{IO}
C. Ionescu, A. Lopes, D. Copot, J.A.T. Machado, J.H.T Bates, The role of fractional calculus in modeling biological phenomena: a review, Commun Nonlinear Sci Numer Simul. 51 (2017) 141--159.


\bibitem{KIL}
A.A. Kilbas, H.M. Srivastava, J.J. Trujillo, Theory and Applications of Fractional Differential Equations, vol. 204, Elsevier Science B.V., Amsterdam, 2006.


\bibitem{KI-NA}
M. Kirane, A. Nabti, Life span of solutions to a nonlocal in time nonlinear fractional 
Schr\"odinger equation, Z. Angew. Math. Phys. 66 (2015) 1473--1482.




\bibitem{LAS}
N. Laskin, Fractional quantum mechanics and L\'evy integral, Phys. Lett. 268 (2000) 298--305.

\bibitem{LAS2}
N. Laskin, Fractional Schr\"odinger equation, Phys Rev E. 66 (2002) 056108.

\bibitem{LI}
M. Li, C. Huang, P. Wang, Galerkin finite element method for nonlinear fractional Schr\"odinger equations, Numer. Algor. 74 (2017) 499--525.

\bibitem{MA}
R.L. Magin, M. Ovadia, Modeling the cardiac tissue electrode interface using fractional calculus, Journal of Vibration and Control. 14 (2008) 1431--1442.

\bibitem{MAI}
F. Mainardi, Fractional relaxation-oscilation and fractional diffusion-wave phenomena, Chaos, Solitons \& Fractals. 7 (1996) 1461--1477.

\bibitem{MI}
E. Mitidieri, S.I. Pohozaev, A priori estimates and blow-up of solutions to nonlinear partial differential equations and inequalities, Proc. Steklov Inst. Math. 234 (2001) 1--383.

\bibitem{NAB}
M.G. Naber, Time fractional Schr\"odinger equation, J. Math. Phys. 45 (8) (2004) 3339--3352.

\bibitem{RA}
M. Ran, C. Zhang, A conservative difference scheme for solving the strongly coupled nonlinear fractional Schr\"odinger equations, Commun. Nonlinear Sci. Numer. 41 (2016)  64--83.

\bibitem{WA}
S. Wang, M. Xu, Generalized fractional Schr\"odinger equation with space-time fractional derivatives, J Math Phys. 48 (2007) 043502.

\bibitem{WA2}
J.R. Wang, Y. Zhou, W. Wei, Fractional Schr\"odinger equations with potential and optimal controls, Nonlinear Anal.: RWA 13 (2012) 2755--2766.

\bibitem{ZH}
Q. Zhang, H.R. Sun, Y. Li, The nonexistence of global solutions for a time fractional nonlinear Schr\"odinger equation without gauge invariance, Appl. Math. Lett. 64 (2017) 119--124.

\end{thebibliography}
\end{document}